

\magnification 1200
\parindent=20pt
\widowpenalty=10000
\def\standardlineskip{12.0pt}

\def\skipfromtitle{\vskip .8cm}
\def\skipfromabstract{\vskip .4cm}

%

\catcode`\@=11

\font\tenmsa=msam10
\font\sevenmsa=msam7
\font\fivemsa=msam5
\font\tenmsb=msbm10
\font\sevenmsb=msbm7
\font\fivemsb=msbm5
\newfam\msafam
\newfam\msbfam
\textfont\msafam=\tenmsa  \scriptfont\msafam=\sevenmsa
  \scriptscriptfont\msafam=\fivemsa
\textfont\msbfam=\tenmsb  \scriptfont\msbfam=\sevenmsb
  \scriptscriptfont\msbfam=\fivemsb

\font\zwartten=msbm7 scaled\magstep2
\font\zwartseven=msbm5 scaled\magstep2
\font\zwartfive=msbm5
\newfam\mszfam
\textfont\mszfam=\zwartten \scriptfont\mszfam=\zwartseven
 \scriptscriptfont\mszfam=\zwartfive
\def\Sb{\ifmmode\let\next\Sb@\else
 \def\next{\errmessage{Use \string\Sb\space only in math mode}}\fi\next}
\def\Sb@#1{{\Sb@@{#1}}}
\def\Sb@@#1{\fam\mszfam#1}

\def\hexnumber@#1{\ifcase#1 0\or1\or2\or3\or4\or5\or6\or7\or8\or9\or
	A\or B\or C\or D\or E\or F\fi }

\font\teneuf=eufm10
\font\seveneuf=eufm7
\font\fiveeuf=eufm5
\newfam\euffam
\textfont\euffam=\teneuf
\scriptfont\euffam=\seveneuf
\scriptscriptfont\euffam=\fiveeuf
\def\frak{\ifmmode\let\next\frak@\else
 \def\next{\errmessage{Use \string\frak\space only in math mode}}\fi\next}
\def\goth{\ifmmode\let\next\frak@\else
 \def\next{\errmessage{Use \string\goth\space only in math mode}}\fi\next}
\def\frak@#1{{\frak@@{#1}}}
\def\frak@@#1{\fam\euffam#1}

\edef\msa@{\hexnumber@\msafam}
\edef\msb@{\hexnumber@\msbfam}

\mathchardef\boxdot="2\msa@00
\mathchardef\boxplus="2\msa@01
\mathchardef\boxtimes="2\msa@02
\mathchardef\square="0\msa@03
\mathchardef\blacksquare="0\msa@04
\mathchardef\centerdot="2\msa@05
\mathchardef\lozenge="0\msa@06
\mathchardef\blacklozenge="0\msa@07
\mathchardef\circlearrowright="3\msa@08
\mathchardef\circlearrowleft="3\msa@09
\mathchardef\rightleftharpoons="3\msa@0A
\mathchardef\leftrightharpoons="3\msa@0B
\mathchardef\boxminus="2\msa@0C
\mathchardef\Vdash="3\msa@0D
\mathchardef\Vvdash="3\msa@0E
\mathchardef\vDash="3\msa@0F
\mathchardef\twoheadrightarrow="3\msa@10
\mathchardef\twoheadleftarrow="3\msa@11
\mathchardef\leftleftarrows="3\msa@12
\mathchardef\rightrightarrows="3\msa@13
\mathchardef\upuparrows="3\msa@14
\mathchardef\downdownarrows="3\msa@15
\mathchardef\upharpoonright="3\msa@16

\mathchardef\downharpoonright="3\msa@17
\mathchardef\upharpoonleft="3\msa@18
\mathchardef\downharpoonleft="3\msa@19
\mathchardef\rightarrowtail="3\msa@1A
\mathchardef\leftarrowtail="3\msa@1B
\mathchardef\leftrightarrows="3\msa@1C
\mathchardef\rightleftarrows="3\msa@1D
\mathchardef\Lsh="3\msa@1E
\mathchardef\Rsh="3\msa@1F
\mathchardef\rightsquigarrow="3\msa@20
\mathchardef\leftrightsquigarrow="3\msa@21
\mathchardef\looparrowleft="3\msa@22
\mathchardef\looparrowright="3\msa@23
\mathchardef\circeq="3\msa@24
\mathchardef\succsim="3\msa@25
\mathchardef\gtrsim="3\msa@26
\mathchardef\gtrapprox="3\msa@27
\mathchardef\multimap="3\msa@28
\mathchardef\therefore="3\msa@29
\mathchardef\because="3\msa@2A
\mathchardef\doteqdot="3\msa@2B

\mathchardef\triangleq="3\msa@2C
\mathchardef\precsim="3\msa@2D
\mathchardef\lesssim="3\msa@2E
\mathchardef\lessapprox="3\msa@2F
\mathchardef\eqslantless="3\msa@30
\mathchardef\eqslantgtr="3\msa@31
\mathchardef\curlyeqprec="3\msa@32
\mathchardef\curlyeqsucc="3\msa@33
\mathchardef\preccurlyeq="3\msa@34
\mathchardef\leqq="3\msa@35
\mathchardef\leqslant="3\msa@36
\mathchardef\lessgtr="3\msa@37
\mathchardef\backprime="0\msa@38
\mathchardef\risingdotseq="3\msa@3A
\mathchardef\fallingdotseq="3\msa@3B
\mathchardef\succcurlyeq="3\msa@3C
\mathchardef\geqq="3\msa@3D
\mathchardef\geqslant="3\msa@3E
\mathchardef\gtrless="3\msa@3F
\mathchardef\sqsubset="3\msa@40
\mathchardef\sqsupset="3\msa@41
\mathchardef\vartriangleright="3\msa@42
\mathchardef\vartriangleleft="3\msa@43
\mathchardef\trianglerighteq="3\msa@44
\mathchardef\trianglelefteq="3\msa@45
\mathchardef\bigstar="0\msa@46
\mathchardef\between="3\msa@47
\mathchardef\blacktriangledown="0\msa@48
\mathchardef\blacktriangleright="3\msa@49
\mathchardef\blacktriangleleft="3\msa@4A
\mathchardef\vartriangle="0\msa@4D
\mathchardef\blacktriangle="0\msa@4E
\mathchardef\triangledown="0\msa@4F
\mathchardef\eqcirc="3\msa@50
\mathchardef\lesseqgtr="3\msa@51
\mathchardef\gtreqless="3\msa@52
\mathchardef\lesseqqgtr="3\msa@53
\mathchardef\gtreqqless="3\msa@54
\mathchardef\Rrightarrow="3\msa@56
\mathchardef\Lleftarrow="3\msa@57
\mathchardef\veebar="2\msa@59
\mathchardef\barwedge="2\msa@5A
\mathchardef\doublebarwedge="2\msa@5B
\mathchardef\angle="0\msa@5C
\mathchardef\measuredangle="0\msa@5D
\mathchardef\sphericalangle="0\msa@5E
\mathchardef\varpropto="3\msa@5F
\mathchardef\smallsmile="3\msa@60
\mathchardef\smallfrown="3\msa@61
\mathchardef\Subset="3\msa@62
\mathchardef\Supset="3\msa@63
\mathchardef\Cup="2\msa@64

\mathchardef\Cap="2\msa@65

\mathchardef\curlywedge="2\msa@66
\mathchardef\curlyvee="2\msa@67
\mathchardef\leftthreetimes="2\msa@68
\mathchardef\rightthreetimes="2\msa@69
\mathchardef\subseteqq="3\msa@6A
\mathchardef\supseteqq="3\msa@6B
\mathchardef\bumpeq="3\msa@6C
\mathchardef\Bumpeq="3\msa@6D
\mathchardef\lll="3\msa@6E

\mathchardef\ggg="3\msa@6F

\mathchardef\circledS="0\msa@73
\mathchardef\pitchfork="3\msa@74
\mathchardef\dotplus="2\msa@75
\mathchardef\backsim="3\msa@76
\mathchardef\backsimeq="3\msa@77
\mathchardef\complement="0\msa@7B
\mathchardef\intercal="2\msa@7C
\mathchardef\circledcirc="2\msa@7D
\mathchardef\circledast="2\msa@7E
\mathchardef\circleddash="2\msa@7F
\def\ulcorner{\delimiter"4\msa@70\msa@70 }
\def\urcorner{\delimiter"5\msa@71\msa@71 }
\def\llcorner{\delimiter"4\msa@78\msa@78 }
\def\lrcorner{\delimiter"5\msa@79\msa@79 }
\def\yen{\mathhexbox\msa@55 }
\def\checkmark{\mathhexbox\msa@58 }
\def\circledR{\mathhexbox\msa@72 }
\def\maltese{\mathhexbox\msa@7A }

\mathchardef\lvertneqq="3\msb@00
\mathchardef\gvertneqq="3\msb@01
\mathchardef\nleq="3\msb@02
\mathchardef\ngeq="3\msb@03
\mathchardef\nless="3\msb@04
\mathchardef\ngtr="3\msb@05
\mathchardef\nprec="3\msb@06
\mathchardef\nsucc="3\msb@07
\mathchardef\lneqq="3\msb@08
\mathchardef\gneqq="3\msb@09
\mathchardef\nleqslant="3\msb@0A
\mathchardef\ngeqslant="3\msb@0B
\mathchardef\lneq="3\msb@0C
\mathchardef\gneq="3\msb@0D
\mathchardef\npreceq="3\msb@0E
\mathchardef\nsucceq="3\msb@0F
\mathchardef\precnsim="3\msb@10
\mathchardef\succnsim="3\msb@11
\mathchardef\lnsim="3\msb@12
\mathchardef\gnsim="3\msb@13
\mathchardef\nleqq="3\msb@14
\mathchardef\ngeqq="3\msb@15
\mathchardef\precneqq="3\msb@16
\mathchardef\succneqq="3\msb@17
\mathchardef\precnapprox="3\msb@18
\mathchardef\succnapprox="3\msb@19
\mathchardef\lnapprox="3\msb@1A
\mathchardef\gnapprox="3\msb@1B
\mathchardef\nsim="3\msb@1C
\mathchardef\ncong="3\msb@1D

\mathchardef\varsubsetneq="3\msb@20
\mathchardef\varsupsetneq="3\msb@21
\mathchardef\nsubseteqq="3\msb@22
\mathchardef\nsupseteqq="3\msb@23
\mathchardef\subsetneqq="3\msb@24
\mathchardef\supsetneqq="3\msb@25
\mathchardef\varsubsetneqq="3\msb@26
\mathchardef\varsupsetneqq="3\msb@27
\mathchardef\subsetneq="3\msb@28
\mathchardef\supsetneq="3\msb@29
\mathchardef\nsubseteq="3\msb@2A
\mathchardef\nsupseteq="3\msb@2B
\mathchardef\nparallel="3\msb@2C
\mathchardef\nmid="3\msb@2D
\mathchardef\nshortmid="3\msb@2E
\mathchardef\nshortparallel="3\msb@2F
\mathchardef\nvdash="3\msb@30
\mathchardef\nVdash="3\msb@31
\mathchardef\nvDash="3\msb@32
\mathchardef\nVDash="3\msb@33
\mathchardef\ntrianglerighteq="3\msb@34
\mathchardef\ntrianglelefteq="3\msb@35
\mathchardef\ntriangleleft="3\msb@36
\mathchardef\ntriangleright="3\msb@37
\mathchardef\nleftarrow="3\msb@38
\mathchardef\nrightarrow="3\msb@39
\mathchardef\nLeftarrow="3\msb@3A
\mathchardef\nRightarrow="3\msb@3B
\mathchardef\nLeftrightarrow="3\msb@3C
\mathchardef\nleftrightarrow="3\msb@3D
\mathchardef\divideontimes="2\msb@3E
\mathchardef\varnothing="0\msb@3F
\mathchardef\nexists="0\msb@40
\mathchardef\mho="0\msb@66
\mathchardef\eth="0\msb@67
\mathchardef\eqsim="3\msb@68
\mathchardef\beth="0\msb@69
\mathchardef\gimel="0\msb@6A
\mathchardef\daleth="0\msb@6B
\mathchardef\lessdot="3\msb@6C
\mathchardef\gtrdot="3\msb@6D
\mathchardef\ltimes="2\msb@6E
\mathchardef\rtimes="2\msb@6F
\mathchardef\shortmid="3\msb@70
\mathchardef\shortparallel="3\msb@71
\mathchardef\smallsetminus="2\msb@72
\mathchardef\thicksim="3\msb@73
\mathchardef\thickapprox="3\msb@74
\mathchardef\approxeq="3\msb@75
\mathchardef\succapprox="3\msb@76
\mathchardef\precapprox="3\msb@77
\mathchardef\curvearrowleft="3\msb@78
\mathchardef\curvearrowright="3\msb@79
\mathchardef\digamma="0\msb@7A
\mathchardef\varkappa="0\msb@7B
\mathchardef\hslash="0\msb@7D
\mathchardef\hbar="0\msb@7E
\mathchardef\backepsilon="3\msb@7F
\def\Bbb{\ifmmode\let\next\Bbb@\else
 \def\next{\errmessage{Use \string\Bbb\space only in math mode}}\fi\next}
\def\Bbb@#1{{\Bbb@@{#1}}}
\def\Bbb@@#1{\fam\msbfam#1}

\catcode`\@=12

\font\fabs=cmr8
\font\bfabs=cmbx8

\font\chaptertitlefont=cmr10 scaled\magstep2
\font\helponefont=cmbx10 scaled\magstep2
\font\helptwofont=cmmi7 scaled\magstep2

\def\pmb#1{\setbox0=\hbox{#1}\copy0\kern-\wd0%
\kern0.02em\raise0.02ex\copy0\kern-\wd0\kern0.02em\box0}

\hyphenation{Hec-ke}

\def\abstract#1{\midinsert\narrower{\noindent\bfabs Abstract. }{\fabs #1}
\endinsert\baselineskip=\standardlineskip}
\def\newsection #1 {\par\noindent{\bf {#1}.}$\,\,\,$}
\def\pr {\noindent{\it Proof}.$\,\,\,$}
\def\Box{\leavevmode\hbox{\vrule\vbox{\hsize=0.6em\hrule
	 \hrule height0.4em depth0.2em width 0pt
	 \hbox{\hskip0.6em plus 0pt minus 1fil}\hrule}\vrule}}
\def\qed{{\unskip\nobreak\hfil\penalty50
	 \hskip2em\hbox{}\nobreak\hfil\Box
	 \parfillskip=0pt \finalhyphendemerits=0 \par
	 \medbreak\noindent\ignorespaces}}

\def\ex #1 {e^{#1}}
\def\st {\,\vert\,}

\def\a {\alpha}
\def\b {\beta}
\def\g {\gamma}
\def\d {\delta}
\def\e {\epsilon}

\def\l {\lambda}
\def\m {\mu}
\def\n {\nu}

\def\p {\pi}

\def\s {\sigma}

\def\ps {\psi}
\def\o {\omega}

\def\S {\Sigma}

\def\O {\Omega}

\def\C {{\bf C}}

\def\R {{\bf R}}

\def\edual {E^\prime}
\def\fialg {{\cal P}_\Phi}
\def\dfialg {{\cal P}_{D\Phi,m}}

\def\schw {{\cal S}}

\def\fourier {{\cal F}}

\def\hol {{\cal H}}

\def\abref {AB}
\def\bref {B}
\def\hamref {H}
\def\rurefone {R1}
\def\rureftwo {R2}
\def\rurefthree {R3}
\def\schwref {S}
\def\treref {T}
\def\zaref {Z}

\centerline{\chaptertitlefont On dense subspaces in a class of}
\centerline{\chaptertitlefont Fr\'echet function spaces on ${\hbox{\helponefont
R}}^{\hbox{\helptwofont n}}$}

\vskip 1cm

\centerline{M.F.E.\ de Jeu}
\centerline{University of Leiden}
\centerline{Department of Mathematics and Computer Science}
\centerline{P.O.\ Box 9512}
\centerline{2300 RA Leiden}
\centerline{The Netherlands}
\centerline{e-mail: jeu@wi.leidenuniv.nl}

\skipfromtitle

\abstract{When dealing with concrete problems in a function space on
${\scriptstyle\R^n}$, it is sometimes helpful to have a dense subspace
consisting of
functions of a particular type, adapted to the problem under
consideration. We give a theorem that allows one to write down many of
such subspaces in commonly occurring Fr\'echet function spaces. These subspaces
are all of the form ${\scriptstyle\{pf_0\st p\in{\cal P}\}}$ where
${\scriptstyle f_0}$ is a fixed function
and ${\scriptstyle{\cal P}}$ is an algebra of functions.
Classical results like the Stone-Weierstrass theorem for polynomials and
the completeness of the Hermite functions are related by this
theorem.}

\skipfromabstract

\xdef\fourone{\the\pageno}

\newsection{1. Introduction}

\vskip .4cm

\par\noindent
Dense subspaces are an important tool in analysis. There are
general theorems guaranteeing the density of subspaces consisting of functions
possessing the greatest possible regularity, such as the density
of $C_c^\infty(\R^n)$ in $L_p(\R^n,dx)\,\,(1\leq p<\infty)$ and in the Schwartz
space
$\schw(\R^n)$. In spite of the significance of these
theorems in the general theory, they are sometimes of little help in more
concrete situations. As an example,
when studying the harmonic oscillator on the real line, one wants to know that
the Hermite functions span a dense subspace of $L_2(\R,dx)$. This is of course
a classical result, but it does not follow from the theorem we mentioned; it
requires a separate proof. We exhibit a theorem that allows a
quick conclusion that some \lq\lq special" subspaces are dense. To illustrate
the idea of the proof, let us sketch it in a particular case.

\vskip .4cm

\par\noindent
Let $\schw(\R)$ be the space of rapidly decreasing functions on the real line,
endowed with its usual Fr\'echet topology. Let $\ps (x)=\exp(-x^2/2)$ and
consider
the subspace $L=\{P\ps\st P$ a polynomial$\}$ of $\schw(\R)$. Then it is known
that $L$ is dense in $\schw(\R)$; this follows e.g.\ from the results in
[\schwref,
p. 263].
The results in [loc.cit.] are based on recurrence relations for the
Hermite functions; we give a more intuitive proof. To this end, fix $T\in
\schw(\R)^\prime$ and consider the map $H_T:\,\C\mapsto\C$, defined by $H_T
(\l)=\langle T\,,\,e_{-i\l}\ps\rangle$ where $e_{-i\l}(x)=\exp(-i\l x)$. Then
$H_T$ is in fact
holomorphic and $({d\over d\l})^k H_T\,(0)=\langle T\,,\,(-ix)^k
\ps\rangle\,\,(k=0,1,\dots)$. Thus, if
$T$ vanishes on $L$, then $\langle T\,,\,e_{-i\l}\ps\rangle=0$ for all
$\l\in\C$, in particular
for all $\l\in\R$. Now due to the completeness of $\schw(\R)$, the weak
integral
$$
I_f=\int_\R f(\l) e_{-i\l}\ps\,d\l
$$
exists in $\schw(\R)$ for all $f\in\schw(\R)$, and is in fact equal to
$\fourier(f) \ps$
(where $\fourier$ denotes Fourier transform). Hence
$$
\langle T\,,\,\fourier(f)\ps\rangle=0\quad(\forall f\in\schw(\R)).
$$
But the Fourier transform maps $\schw(\R)$ onto itself, so
$$
\langle T\,,\,f\ps\rangle=0\quad(\forall f\in\schw(\R)).
$$
Now observe that $\{f\ps\,|\,f\in\schw(\R)\}$ is dense in $\schw(\R)$: it
contains
$C_c^\infty(\R)$ since $\ps$ has no zeros. We conclude that
$T=0$ and, finally, that $L$ is dense in $\schw(\R)$ by the Hahn-Banach
theorem.

\vskip .4cm

\par\noindent
The above proof is based on a combination of function theory and
Fourier analysis. The application of this combination in density problems
has a long history: it goes at
least back to Hamburger's work in 1919 on $L_2((0,\infty),dx)$ ([\hamref]).
This paper
fits into this tradition: it turns out that the combination of function theory
and Fourier analysis
can be put to good use in a more general context to supply dense subspaces
in function spaces on $\R^n$, provided that the topology is defined in a
certain way (to be explained in Section 2). We (must) assume that the space is
Fr\'echet, since the existence of vector-valued integrals as in the example
above is essential.

\vskip .4cm
\par\noindent
The main theorem is Theorem 2.13, stating that the annihilators of certain
subspaces are equal. In particular,
one of them is dense if and only if the other is. E.g., in the above example
it follows from Theorem 2.13 that $\{P\ps\st P$ a polynomial$\}$ and $\{f\ps\st
f\in\schw(\R)\}$ have the same annihilator, and we happen to know that the
latter
subspace is dense in $\schw(\R)$ since it contains $C_c^\infty(\R)$. Thus it is
the
combination of Theorem 2.13 below and \lq\lq general" density theorems that
allows one
to conclude that some \lq\lq special" subspaces are dense as well.

\vskip .4cm

\par\noindent
This paper is organized as follows. In Section 2, we start with the observation
that the way in which a number of well-known function spaces are topologized
can be described in a uniform manner. This being done, we prove the main
theorem. Section 3 contains applications in three cases. These cases
do not exhaust the possible applications of the method in
this paper; they rather serve as an illustration, leaving it to the reader to
apply the method to situations of his interest. We conclude
in Section 4 with remarks on possible variations of the method and connections
with
representation theory.

\vskip .6cm

\xdef\fourtwo{\the\pageno}
\newsection{2. Main theorem}

\vskip .4cm

\par\noindent
The proof in the Introduction works for a whole class of function spaces
on $\R^n$, provided that the topology is defined in a certain manner. After
establishing the conventions and notation, we start with an attempt to
formalize the way in
which a number of function spaces are topologized, and we give some examples.
We then work towards Theorem 2.13 on equality of annihilators.

\vskip .4cm

\par\noindent
All topological function spaces under consideration are assumed to be complex.
This is essential during the proof (since holomorphic functions are involved),
but in applications it is usually an easy matter to derive a result for the
real case from the result for the complex case. In order to guarantee the
existence
of certain vector-valued integrals, we assume that the spaces are Fr\'echet
(which includes local convexity by convention) from the start, although the
results up to and including Corollary 2.11 also hold without this assumption.
\par\noindent
By convention, Borel measures take finite values on compact sets, which implies
that any extension of a Borel measure on an open subset of $\R^n$ is
$\s$-finite (a technical condition that allows application of
Fubini's theorem).
\par\noindent
The argument of functions is usually omitted;
(in)equalities involving functions should always be read pointwise almost
everywhere.

\vskip .4cm

\par\noindent
We write $(.\,,\,.)$ for the usual bilinear form on $\C^n\times\C^n$: if
$x=(x_1,\dots,x_n)$ and $y=(y_1,\dots,y_n)$ are in $\C^n$, then
$(x,y)=\sum_{j=1}^nx_jy_j$. The standard two-norm on $\C^n$ is denoted by
$\Vert\,.\,\Vert$; we have $\vert(x,y)\vert\leq\Vert x\Vert\,\Vert y\Vert$.
\par\noindent
The standard multi-index notation is used throughout: if
$\a=(\a_1,\dots,\a_n)$, then $\vert\a\vert=\sum_{i=1}^n\a_i$ and
$D^\a=\left({\partial\over\partial x_1}\right)^{\a_1}\dots\left(
{\partial\over\partial x_n}\right)^{\a_n}
$. We write $\b\leq\a$ if $\b_i\leq\a_i\,\,(i=1,\dots\,n)$.
\par\noindent
For $\e>0$, let $\C_\e^n=\{\l\in\C^n\st \Vert{\rm Im}\l\Vert<\e\}$.
\par\noindent
If $E$ is a topological vector space, then $E^\prime$ denotes its dual. For
$S\subset E$ we let $S^\perp=\{T\in E^\prime\st Ts=0\,\,\forall s\in S\}$ be
the annihilator of $S$.
\par\noindent
The Fourier transform, finally, is denoted by $\fourier$.

\vskip .4cm

\par\noindent
We observe that the topology of many locally convex function spaces on $\R^n$
is defined by seminorms involving only one or more of the following
ingredients:
\item{--} amount of differentiability,
\item{--} amount of integrability with respect to a measure,
\item{--} behaviour on compacta, and
\item{--} weight functions.
\par\noindent
Moreover, many of these spaces are in fact metrizable and complete.
The following definition introduces an ad hoc terminology for this kind of
spaces and formalizes a common method of topologizing function spaces on
$\R^n$.

\vskip .4cm

\par\noindent
\newsection {Definition 2.1} A complex Fr\'echet space $E$ is a {\it common
Fr\'echet
function space on $\R^n$} if there is a sextuple $(U,\{U_k\}_{k=1}^\infty,
\m,p,m,\{\nabla_N^\a\}_{N=0,1,2,\dots\,;\vert\a\vert\leq m})$ such that:
\vskip .2cm
\item{1.} $U\subset\R^n$ is open and non-empty.
\item{2.} $U_k\subset\R^n$ is open for all $k$, and $U=\bigcup_{k=1}^\infty
U_k$.
\item{3.} $\m$ is the completion of a Borel measure on $U$. Let $M(U)$
          denote the vector space of $\m$-measurable functions on $U$, where we
agree to
          identify two elements if they are equal a.e. $(\m)$.
\item{4.} $E$ is a subspace of $M(U)$.
\item{5.} $1\leq p\leq\infty$.
\item{6.} $m\in\{0,1,2,\dots\}\cup\{\infty\}$.
\item{7.} The $\nabla_N^\a$ are linear maps $\nabla_N^\a:E\mapsto M(U)$,
          indexed by a non-negative integer $N$ and a multi-index $\a=(\a_1,
          \dots,\a_n)$ subject to the condition $\vert\a\vert\leq m$, such
that:
          \itemitem{--}$\nabla_0^0$ is the inclusion of $E$ in $ M(U)$;
          \itemitem{--}$
          \Vert\chi_k\nabla_N^\a e\Vert_p<\infty
          \quad\left(\forall e\in E,\,\,\forall k,\,\,\forall N,\,\,\forall
\a\,\,(\vert\a\vert\leq m)\right).
          $
          Here $\chi_k$ denotes the characteristic function of $U_k$ and the
          $p$-norm corresponds to the measure $\m$.
\item{8.} The topology on $E$ is defined by seminorms
          $\{p_{k,\a,N}\}_{k,N=0,1,2,\dots\,;\vert\a\vert\leq m}$ on $E$, where
          $p_{k,\a,N}$ is defined by $p_{k,\a,N}(e)
          =\Vert\chi_k\nabla_N^\a e\Vert_p\,\,(e\in E)$.
\item{9.} For all $\a\,\,(\vert\a\vert\leq m)$ and $\b,\g\leq\a$ there exist
constants
          $c_{\b,\g;\a}$ such that
          $$
          \nabla_N^\a(g e)=\sum_{\b+\g=\a}c_{\b,\g;\a}\,
          (D^\b g)(\nabla_N^{\g}e)
          \quad (\forall N)
          $$
          whenever $g\in C^\infty (U)$ and $e\in E$ are such that $ge\in E$.

\vskip .4cm
\par\noindent
The integer $m$ should be thought of as describing that the distributional
derivatives up to order
$m$ are \lq\lq regular" (this usually means that the the function is $C^m$, but
it has
a different meaning in the context of Sobolev spaces), the $\nabla_N^\a$ can
be interpreted as distributional differentiation followed by multiplication by
a weight
function, the $U_k$ allow incorporation of behaviour on compacta, and $p$ of
course expresses integrability (essential boundedness if $p=\infty$). The
following
examples illustrate this.

\vskip .4cm

\par\noindent
\newsection{Example 2.2} Let $\m$ be the completion of a Borel measure on
an open subset $U$ of $\R^n$. Then $L_p(U,\m)\,\,(1\leq p\leq\infty)$
is a common Fr\'echet function space: put $U_k=U\,\,(\forall k),\,m=0,\,$ and
let $\nabla_N^0$ be the inclusion for all $N$.

\vskip .4cm
\newsection{Example 2.3} Let $U\subset\R^n$ be open. The space $C^m(U)\,\,
(1\leq m\leq\infty)$ is
canonically embedded in $M(U)$ (here Lebesgue measure is tacitly understood).
Choose a sequence $\{U_k\}_{k=1}^\infty$ of open subsets of $U$ such that
$U=\bigcup_{k=1}^\infty U_k$, $\overline{U_k}$ is compact for all $k$, and
$\overline{U_k}\subset U_{k+1}$ for all $k$. Put
$\nabla_N^\a=D^\a\,\,(\vert\a\vert
\leq m,\,N\geq 0)$ (note that this a legitimate definition on the
embedding of $C^m(U)$, since an equivalence class contains exactly one
continuous representative). Let $p=\infty$. Then the usual topology on
$C^m(U)$ is obtained, showing that $C^m(U)$ is a common Fr\'echet function
space.

\vskip .4cm

\newsection{Example 2.4} Let $\schw(\R^n)$ be the space of all smooth functions
of
rapid decrease at infinity, canonically embedded in $M(\R^n)$ (where Lebesgue
measure
is again understood). Let
$U_k=\R^n\,\,(\forall k)$, and put $(\nabla^\a_N f)(x)=(1+\Vert x\Vert)^N D^\a
f(x)
\,\,(x\in\R^n,\,f\in\schw(\R^n))$. Let $p=m=\infty$. Thus $\schw(\R^n)$ is a
common Fr\'echet
function space.

\vskip .4cm

\par\noindent
The reader will have no trouble verifying that other spaces (e.g.,
$L_p^{loc}(U)$
and the Sobolev spaces $W^{m,p}(U)$) are common in the sense of the above
definition.

\vskip .4cm

\par\noindent
We now embark on the proof of the main result, Theorem 2.13. The proof consists
of two steps, as in the Introduction. The first step consists of showing
that certain $E$-valued maps are holomorphic, and in the second one an
$E$-valued
integral is identified. The two results together then easily yield the theorem.
\vskip .4cm
\par\noindent
As it turns out, the polynomials in the example in the Introduction can be
replaced by polynomial functions of the components of an arbitrary smooth map
$\Phi:\,U\mapsto\R^n$ (a diffeomorphism in the applications in Section 3)
without complicating the proof. Since derivatives up to order $m$ also appear
in the topology of $E$ (notably in Definition 2.1.9),
we are led to the following definition.

\vskip .4cm
\par\noindent
\newsection{Definition 2.5} Let $\Phi:\,U\mapsto\R^n$ be of class $C^\infty$
with components $\Phi_1,\dots,\Phi_n$, let $\e>0$ and let $m$ be a non-negative
integer.
Then:
\item{1.} $\fialg$ is the unital algebra of functions on $U$ generated by
          $\{\Phi_j\st 1\leq j\leq n\}$.
\item{2.} $\dfialg$ is the unital algebra of functions on $U$ generated by
         $\{D^\a\Phi_j\st 1\leq j\leq n,\,\,\vert\a\vert\leq m \}$.

\vskip .4cm

\par\noindent
In order to avoid repetitions, we make the following assumption for the
remainder of this section.

\vskip .4cm

\par\noindent
\newsection{Assumption 2.6}  We assume that we are given a common complex
Fr\'echet
function space $E$ with defining
sextuple $(U,\{U_k\}_{k=1}^\infty,\m,p,m,\{\nabla_N^\a\}_{N=0,1,2,\dots\,
;\vert\a\vert\leq m})$, a $C^\infty$-map
$\Phi:\,U\mapsto\R^n$, $e_0\in E$, and $\e>0$ such that:
\item{1.} $\ex {i(\l,\Phi)} ge_0\in E\,\,(\forall\l\in\C^n_\e,\,\forall
          g\in\fialg)$.
\item{2.} $\left\Vert\chi_k \ex {\e\Vert\Phi\Vert} g\nabla_N^\a
e_0\right\Vert_p
          <\infty\quad(\forall k,\forall\a\,(\vert\a\vert\leq m),
          \forall N,\,\forall g\in\dfialg)$.

\vskip .4cm
\par\noindent
To obtain the example in the Introduction, one takes $E=\schw(\R^n)$, $\Phi$
the identity,
$e_0(x)=\exp (-x^2/2)$ and any $\e>0$.

\vskip .4cm

\par\noindent
For $p\in \fialg$, define the map $\hol_{p}:\C_\e^n\mapsto E$
by:
$$
\hol_{p}(\l)=\ex {-i(\l,\Phi)} p e_0\quad(\l\in\C_\e^n).
$$
We start by proving that $\hol_{1}$ is weakly holomorphic and identifying
its derivatives (Proposition 2.10).

\vskip .4cm

\par\noindent
\proclaim Lemma 2.7. $\hol_{p}$ is continuous on $\C_\e^n$ for all
$p\in\fialg$.

\vskip .4cm

\pr
Fix $\l^0\in\C_\e^n$. We are to prove that
$$
\lim_{h\rightarrow 0}p_{k,N,\a}\left(\ex {-i(\l^0+h,\Phi)}
p e_0-
\ex {-i(\l^0,\Phi)} p e_0\right)=0\quad(\forall k,\,\,\forall
N,\,\,\forall\a\,\,(\vert\a\vert\leq m)).
$$
Fix $k,N$ and $\a$. We have
$$
\eqalign
{
&p_{k,N,\a}\left(\ex {-i(\l^0+h,\Phi)}
p e_0-
\ex {-i(\l^0,\Phi)} p e_0\right)=\cr
&=\left \Vert\sum_{\b+\g=\a}c_{\b,\g;\a}\,\chi_k\, D^\b\left\{
\ex {-i(\l^0,\Phi)} \left[\ex {-i(h,\Phi)} -1 \right]p\right\}
\nabla_N^{\g}e_0\right\Vert_p.
\cr
}
$$
The summation can be written as a finite sum of terms of two types
(corresponding
to zero and one or more differentiations acting on the expression in square
brackets, respectively):
$$
\chi_k\ex {-i(\l^0,\Phi)} \left(\ex {-i(h,\Phi)} -1 \right)g\nabla_N^\d e_0
\eqno {{\rm (type~1) }}
$$
(where $g\in\dfialg$ does not depend on $h$, and $\vert\d\vert\leq m$), and
$$
\chi_k Q_{g,\d}(h) \ex {-i(\l^0,\Phi)} \ex {-i(h,\Phi)} g\nabla_N^\d e_0
\eqno {{\rm (type~2) }}
$$
(where again $g\in\dfialg$ does not depend on $h$, $\vert\d\vert\leq m$, and
$Q_{g,\d}$ is a polynomial on $\C^n$ such that $Q_{g,\d}(0)=0$).
To estimate terms of the first type, note that
$$
\ex {-iz} -1=-\int_0^1 iz\ex {-itz} dt\quad(z\in\C),\leqno (1)
$$
hence
$$
\vert \ex {-i(h,\Phi)} -1 \vert\leq\Vert h\Vert\,\Vert\Phi\Vert\,
\ex {\Vert{\rm Im}\,h\Vert\,\Vert\Phi\Vert} .
$$
This readily implies that terms of the first type are dominated by
$$
\Vert h \Vert\,\left\vert\chi_k\ex {\e\Vert\Phi\Vert} \Vert\Phi\Vert\,g\,
\nabla_N^\d e_0\right\vert
\leq
\Vert h \Vert \,\sum_{j=1}^n\left\vert\chi_k\,\ex {\e\Vert\Phi\Vert} \Phi_j
g\nabla_N^\d e_0\right\vert
$$
if $\Vert{\rm Im}\,h\Vert + \Vert {\rm Im}\,\l^0\Vert\leq\e$.
\par\noindent
Under this same condition for $h$, terms of the second type are dominated by
$$
\vert Q_{g,\d}(h)\vert \, \left\vert\chi_k \ex {\e\Vert\Phi\Vert} g\,
\nabla_N^\d e_0\right\vert.
$$
Now note that an inequality $\vert f\vert\leq\sum_{j=1}^N \vert f_j\vert$
implies that
$\Vert f\Vert_p\leq\sum_{j=1}^N \Vert f_j\Vert_p$.
Hence the lemma follows from Assumption 2.6.2 and the observation that each of
the majorants contains a
term that tends to zero as $h$ tends to zero.

\qed

\vskip .4cm

\par\noindent
We need the following lemma to establish the existence of $E$-valued integrals
later on. The easy proof is left to the reader.

\vskip .4cm

\par\noindent
\proclaim Lemma 2.8. For all $k,\,N$ and $\a\,\,(\vert\a\vert\leq m)$
there exists a polynomial $Q_{k,N,\a}$ in one variable such that
$$
p_{k,N,\a}(\hol_{1}(\l))\leq Q_{k,N,\a}(\Vert\l\Vert)
\quad (\forall\l\in\C_\e^n).
$$

\vskip .4cm

\proclaim Lemma 2.9. For all $p\in\fialg$, $\hol_{p}$ is holomorphic on
$\C_\e^n$ in each variable
separately and
$$
{d\over d\l_l}\hol_{p}=
\hol_{-i\Phi_l p}
\quad (l=1,\dots,n).
$$

\vskip .4cm

\pr
Let us prove the lemma for the first variable. Let $e_1=(1,0,\dots,0)\in\C^n$;
fix $\l^0\in\C^n_\e ,\,k,N$ and $\a$ and consider for $h\in\C\,\,(h\neq 0)$:
$$
\eqalign
{
&p_{k,N,\a}\left({\ex {-i(\l^0+he_1,\Phi)} p e_0-\ex {-i(\l^0,\Phi)} p
e_0\over h}+i\ex {-i(\l^0,\Phi)} \Phi_1 p e_0 \right)=
\cr
&=\left\Vert\sum_{\b+\g=\a}c_{\b,\g;\a}\chi_k D^\b\left\{\ex {-i(\l^0,\Phi)}
\left[{\ex {-ih\Phi_1} -1\over
h}+i\Phi_1\right]\,p\right\}\nabla_N^{\g}e_0\right\Vert_p.
\cr
}
$$
Now a moment's thought shows that the summation can be written as a finite sum
of terms of three types (corresponding to zero, one, and two or more
differentiations
acting on the expression in square brackets, respectively):
$$
\chi_k \ex {-i(\l^0,\Phi)} \left\{ {\ex {-ih\Phi_1} -1\over
h}+i\Phi_1\right\}\, g
\nabla_N^\d e_0
\eqno {{\rm (type~1) }}
$$
(where $g\in\dfialg$ does not depend on $h$, and $\vert\d\vert\leq m$),
$$
\chi_k \ex {-i(\l^0,\Phi)} \left(\ex {-ih\Phi_1} -1 \right) g
\nabla_N^\d e_0
\eqno {{\rm (type~2) }}
$$
(where again $g\in\dfialg$ does not depend on $h$, and $\vert\d\vert\leq m$),
and
$$
\chi_k Q_{g,\d}(h)\ex {-i(\l^0,\Phi)} \ex {-ih\Phi_1} g\nabla_N^\d e_0
\eqno {{\rm (type~3) }}
$$
(where
$g\in\dfialg$ does not depend on $h$, $\vert\d\vert\leq m$,
and $Q_{g,\d}$ is a polynomial on $\C$ such that $Q_{g,\d}(0)=0$).
\par\noindent
Now note that
$$
\ex {-iz} -1 +iz=-\int_0^1 z^2 t\ex {i(t-1)z} \,dt\quad(z\in\C),
$$
hence
$$
\eqalign
{
\left\vert {\ex {-ih\Phi_1} -1\over h} +i\Phi_1\right\vert&\leq
\vert h\vert\,\vert\Phi_1\vert^2\,\ex {\vert {\rm
Im}\,h\vert\,\vert\Phi_1\vert} \cr
&\leq \vert h\vert\,\vert\Phi_1\vert^2\,\ex {\vert {\rm Im}\,h\vert\,\Vert\Phi
\Vert} .\cr
}
$$
Thus, if $\vert{\rm Im}\,h\vert+\Vert{\rm Im}\,\l^0\Vert\leq\e$, then each
term of the first type is dominated by
$$
\vert h\vert\,
\left\vert \chi_k\,\ex {\e\Vert\Phi\Vert}
\Phi_1^2 g\nabla_N^\d e_0\right\vert.
$$
If one uses (1) again, it is easy to see
that under the same condition for $h$ the terms of the second type are
dominated by
$$
\vert h\vert\,\left\vert\chi_k\ex {\e\Vert\Phi\Vert}
\Phi_1 g\nabla_N^\d e_0\right\vert.
$$
Terms of the third type are dominated by
$$
\vert Q_{g,\d}(h)\vert\,\left\vert\chi_k\,
\ex {\e\Vert\Phi\Vert} \,g\nabla_N^\d e_0\right\vert.
$$
Each of the majorants contains a term that tends to zero if $h$ tends to zero,
so the lemma follows as in the conclusion of the proof of Lemma 2.7.
\qed

\vskip .4cm

\proclaim Proposition 2.10. Under Assumption 2.6 the map
$T\circ\hol_{1}:\C_\e^n\mapsto\C$
is holomorphic for all $T\in E^\prime$, and
$$
D^\a(T\circ\hol_{1})=
T\circ\hol_{(-i)^{\vert\a\vert}\Phi_1^{\a_1}\dots\Phi_n^{\a_n}}
\quad(\forall\a).
$$

\vskip .4cm

\pr
The weak holomorphy of $\hol_{1}$ follows from Lemma 2.7 and
Lemma 2.9 (or from Lemma 2.9 alone if one invokes Hartog's theorem).
The derivatives are identified by repeated application of Lemma 2.9.
\qed

\vskip .4cm

\proclaim Corollary 2.11. Under Assumption 2.6 we have
$$
\eqalign
{
(\fialg e_0)^\perp&=({\rm Span}\{\ex {i(\l,\Phi)}  e_0\st\l\in\C_\e^n \})^\perp
\cr
&=({\rm Span}\{\ex {i(\l,\Phi)}  e_0\st\l\in\R^n \})^\perp.
\cr
}
$$

\vskip .4cm

\par\noindent
The completeness of $E$ (which has not been used until now) is brought into
play in the following lemma.

\vskip .4cm

\proclaim Lemma 2.12. Under Assumption 2.6  the weak
integral
$$
I_f=\int_{\R^n}f(\l)\hol_{1}(\l)\,d\l
$$
exists in $E$ for all $f\in\schw(\R^n)$, and is equal to
$\left(\fourier(f)\circ\Phi\right)
e_0$.

\vskip .4cm

\pr
We recall a basic existence theorem for weak integrals [\rurefone, Theorem 3.27
and
the remark preceding the theorem]: if $E$ is a Fr\'echet space, $X$ is a
compact
Hausdorff space, $\Psi:\,X\mapsto E$ is continuous and $\m$ is a bounded
Borel measure $X$, then the weak integral $\int_X \Psi\,d\m$
exists in $E$. Lemma 2.8 enables one to invoke this theorem, as follows. Let
$X=\R^n\cup\{\infty\}$ be the one-point compactification of $\R^n$. As a
consequence of Lemma 2.8 and the fact that $f\in\schw(\R^n)$, the map $\l
\mapsto (1+\Vert\l\Vert)^N f(\l)\hol_{1}(\l)$ from $\R^n$
into $E$ extends to a continuous map $\Psi_N:\,X\mapsto E$ by putting
$\Psi(\infty)
=0$, for any $N\geq 0$. Choose $N$ so large that $\int_{\R^n} (1+\Vert\l\Vert)
^{-N}\,d\l<\infty$. Extend the measure $(1+\Vert\l\Vert)^{-N}\,d\l$ from $\R^n$
to
a bounded measure on $X$ by declaring the measure of $\{\infty\}$ to be zero.
Now apply the existence theorem to $\Psi_N$.
\par\noindent
To identify the integral, fix $k$ and note that
$\chi_k e\in L_p(U_k,\m)$ for all $e\in E$, since (in the notation of
Definition 2.1) $p_{k,0,0}(e)=\Vert\chi_k e\Vert_p$ is one of the seminorms
defining the topology on $E$. Hence any $g\in
L_q(U_k,\m)$ (where $q$ is the conjugate exponent of $p$) defines
an element of $\edual$. So by the very definition of a weak integral we have
$$
\langle g\,,\,I_f\rangle=\int_{\R^n}f(\l)\,\langle g\,,\,\ex {-i(\l,\Phi)}
e_0\rangle
\,d\l,
$$
where
$$
\langle g\,,\, \ex {-i(\l,\Phi)} e_0\rangle=
\int_{U_k}g(x)\,\ex {-i(\l,\Phi(x))} e_0(x)\,d\m (x).
$$
Now an application of Fubini's theorem shows that the integrals
$$
\int_{\R^n}f(\l)\left\{\int_{U_k}g(x)\,\ex {-i(\l,\Phi(x))} e_0(x)\,d\m(x)
\right\}\,d\l
$$
and
$$
\int_{U_k}g(x)\left\{\int_{\R^n}f(\l)\,\ex {-i(\l,\Phi(x))} \,d\l\right\}e_0(x)
d\m(x)
$$
are equal, i.e.
$$
\langle g\,,\,I_f\rangle=\int_{U_k}g(x)\left(\fourier(f)\circ\Phi\right)(x)
e_0(x)\,d\m(x)\quad (\forall g\in L_q(U_k,\m)).
$$
Now recall that an $L_p$-space is always canonically embedded in $L_q^\prime$
if $1\leq p
<\infty$, and that this embedding also holds for $p=\infty$ if the measure is
$\s$-finite ([\zaref, Lemma $\b$, p. 357]). We conclude that $I_f=
(\fourier(f)\circ\Phi)\,e_0$ a.e. $(\m)$ on $U_k$, which proves the lemma since
the $U_k$ cover $U$.
\qed

\vskip .4cm

\par\noindent
We finally arrive at the theorem on equality of annihilators that was mentioned
in the Introduction.

\vskip .4cm

\proclaim Theorem 2.13. Under Assumption 2.6 the
following subspaces of $\edual$ are equal:
\item{1.}$\{\left(\fourier(f)\circ\Phi\right) e_0\st f\in\schw(\R^n)\}^\perp$.
\item{2.}$\{\left(f\circ\Phi\right) e_0\st f\in\schw(\R^n)\}^\perp$.
\item{3.}$(\fialg e_0)^\perp$.
\item{4.}$({\rm Span}\{\ex {i(\l,\Phi)}  e_0\st\l\in\C_\e^n \})^\perp$.
\item{5.}$({\rm Span}\{\ex {i(\l,\Phi)}  e_0\st\l\in\R^n \})^\perp$.

\vskip .4cm

\pr
The equality of 3, 4 and 5 is just Corollary 2.11.
Since $\fourier:\schw(\R^n)\mapsto\schw(\R^n)$ is a bijection, the equality of
1 and
2 follows trivially. The inclusion 5$\subset$1 is a consequence of Lemma 2.12.
As to the converse, suppose that $T\in\{\left(\fourier(f)\circ\Phi\right)
 e_0\st f\in\schw(\R^n)\}^\perp$. Then Lemma 2.12 shows that
$$
\int_{\R^n}f(\l)(T\circ\hol_{1})(\l)\,d\l=0
\quad(\forall f\in\schw(\R^n)).
$$
Since $T\circ\hol_{1}$ is continuous, it is identically zero,
proving 1$\subset$5.
\qed

\vskip .4cm

\xdef\fourthree{\the\pageno}
\newsection {3. Applications}

\vskip .4cm

\par\noindent
In this section, we obtain some density results for the spaces in the
Examples 2.2-2.4 by combining Theorem 2.13 and the
\lq\lq standard" density theorems. As we remarked in the Introduction, the
applications in this section merely serve to illustrate how
this combination can be used to conclude that subspaces of a
certain type are dense.
We therefore emphasize that we apply
Theorem 2.13 rather crudely by requiring $\Phi:\,U\mapsto
\Phi (U)$ to be a diffeomorphism. Under this condition, the subspace of
pull-backs
$\{f\circ\Phi\st f\in\schw(\R^n)\}$ (figuring in Theorem 2.13) contains
$C_c^\infty(U)$, which is usually \lq\lq large". This enables us to invoke
standard
density theorems. In situations where $\Phi$ is not a diffeomorphism, the
reader
may still have some use for Theorem 2.13, depending on the particular
circumstances
in the problem under consideration.

\vskip .4cm

\par\noindent
As in the previous section, all topological vector space are assumed to
be complex. The reader will have no trouble deriving theorems for the real
case from the results for the complex case.

\vskip .4cm

\par\noindent
We must distinguish between cases at this stage. The reason for this is more
or less obvious: if e.g.\ $E=\schw(\R^n)$, then
the existence of a zero for $e_0$ implies that none of the annihilators in
the Theorem 2.13 is zero, whereas a statement on the existence of a single zero
is generally simply meaningless in the case of $L_p(U,\m)$.
We therefore elaborate separately on the Examples 2.2-2.4.
\par\noindent
For the convenience of the reader, we recall some notation
from the previous section: if $\e>0$, then $\C_\e^n=\{\l\in\C^n\st \Vert{\rm
Im}\,\l
\Vert <\e\}$, and if $\Phi: U\mapsto\R^n$ has components $\Phi_1,\dots,\Phi_n$
then $\fialg$ is the unital algebra of functions on $U$ generated by
$\{\Phi_j\st 1\leq j\leq n\}$.

\vskip .4cm

\par\noindent
{\bf Theorem 3.1.} {\sl Let $\m$ be the completion of a Borel measure on an
open subset $U$ of $\R^n$, and let $1\leq p<\infty$. Let $\Phi:U\mapsto\Phi(U)$
be
a diffeomorphism of class $C^\infty$. Let $f_0\in L_p(U,\m)$ and suppose that
there exists
$\e>0$ such that
$$
\left\Vert \ex {\e\Vert\Phi\Vert} p\,f_0\right\Vert_p<\infty
\quad(\forall p\in\fialg).
$$
Then the annihilators of the following subspaces of $L_p(U,\m)$:
\item{1.}$\{\left(f\circ\Phi\right)f_0\st f\in\schw(\R^n)\}$.
\item{2.}$\fialg f_0$.
\item{3.}${\rm Span}\{\ex {i(\l,\Phi)} f_0\st\l\in\C_\e^n \}$.
\item{4.}${\rm Span}\{\ex {i(\l,\Phi)} f_0\st\l\in\R^n \}$.
\par\noindent
are all equal to $\{g\in L_q(U,\m)\st gf_0=0 \hbox{\rm ~a.e. }(\m)\}$, where
$q$ is the conjugate exponent of $p$.
In particular, these subspaces are dense in $L_p(U,\m)$ if and only if
$f_0(x)\neq 0$
for almost all $x\,\,(\m)$.}

\vskip .4cm
\pr
The equality of the annihilators is just an application of Theorem 2.13.
If $gf_0=0$ a.e. ($\m$), then $g$ is obviously in the annihilator. Conversely,
let $g\in L_q(U,\m)$ be in the annihilator. The subspace
in 1 contains $C_c^\infty(U)f_0$, so in particular
$$
\int_{U} g\,h\,f_0\,d\m=0\quad(\forall h\in\C_c^\infty(U)).
$$
Since $C_c^\infty(U)$ is dense in $C_0(U)$ (the continuous functions on $U$
vanishing at infinity) under the
supremum-norm, the dominated convergence theorem shows that
$$
\int_{U} g\,f_0\,h\,d\m=0\quad(\forall h\in C_0(U)).
$$
Hence $\int_U \vert gf_0 \vert\,d\m=0$ by the Riesz representation theorem
([\rurefthree, Theorem 6.19]), i.e.\ $gf_0=0$ a.e. $(\m)$.
\par\noindent
The criterium for density is a direct consequence of the description
of the annihilator.
\qed

\vskip .4cm

\par\noindent
The theorem also holds if $\m$ is simply a Borel measure, or any measure $\n$
that is intermediate (in the sense of domains of definition and extension)
between a Borel measure and its completion. Indeed, suppose that
$\m$ is a Borel measure on $U$ with the $\s$-algebra ${\cal B}(U)$ of
Borel subsets of $U$ as domain of definition, and let $\m^*$ be its
completion with corresponding $\s$-algebra ${\cal B}^*(U)$. Suppose that
$\S$ is a $\s$-algebra such that ${\cal B}(U)\subset\S\subset{\cal B}^*
(U)$, and let $\n$ be a measure on $\S$ such that the restriction of
$\n$ to ${\cal B}(U)$ is equal to $\m$. Then it follows from [\abref, p. 92]
that $\n$ is in fact equal to the restriction of $\m^*$ to $\S$. But then
[\rurefthree, Lemma 1, p. 154] shows that the natural embedding of $L_p(U,\n)$
in $L_p(U,\m^*)$ is in fact surjective, hence an isometric isomorphism.
Thus, since the theorem holds for $\m^*$, it also holds for $\n$.

\vskip .4cm

\par\noindent
To illustrate the theorem, consider the case $n=1$. Take $f_0(x)=\ex
{{-x^2}\over 2} $, let $\m$ be
Lebesgue measure and let $\Phi$ be the identity. Then Theorem 3.1 asserts that
the Hermite functions span a dense subspace in $L_p(\R,dx)\,\,(1\leq
p<\infty)$. If we restrict the Lebesgue
measure to $(0,\infty)$ and put $f_0(x)=\ex {-x\over 2} x^{\a\over
2}\,\,(\a>-1)$,
then the density of the span of the Laguerre functions in $L_p((0,\infty),dx)$
is obtained
for $1\leq p<\infty$ if $\a\geq 0$, and for $1\leq p<-2/\a$ if $-1<\a<0$. In
particular,
this span is dense in $L_2((0,\infty),dx)$ for all $\a>-1$.
The theorem also allows one to conclude that more exotic subspaces are dense,
e.g.\ (with $[\,.\,]$ denoting the entier function)
$$
\{P(x\sqrt{1+x^2})\ex {-\sqrt {x^6+\cos x + 2}} \,[x^2+2] \st
P{\rm~a~polynomial}\}
$$
is dense in $L_p(\R,dx)$ if $1\leq p<\infty$.

\par\noindent
For arbitrary $n$, the polynomials are dense
in $L_p(U,\m)\,\,(1\leq p<\infty)$ if there exists $\e>0$ such that
$\ex {\e\Vert x\Vert} \in L_p(U,\m)$ (take $f_0=1$).

\vskip .4cm

\par\noindent
{\bf Theorem 3.2.} {\sl Let $U\subset\R^n$ be open and let $0\leq m\leq\infty$.
Endow $C^m(U)$ with its usual topology of uniform convergence of all
derivatives of order $\leq m$ on compact subsets of $U$. Let $f_0\in
C^\infty(U)$, and let $\Phi:U\mapsto\Phi(U)$ be a diffeomorphism of
class $C^\infty$.
Then the annihilators of the following subspaces of $C^m(U)$ are equal:
\item{1.}$\{\left(f\circ\Phi\right)f_0\st f\in\schw(\R^n)\}$.
\item{2.}$\fialg f_0$.
\item{3.}${\rm Span}\{\ex {i(\l,\Phi)} f_0\st\l\in\C_\e^n \}$.
\item{4.}${\rm Span}\{\ex {i(\l,\Phi)} f_0\st\l\in\R^n \}$.
\par\noindent
These subspaces are dense in $C^m(U)$ if and only if $f_0$ has no zeros.}

\vskip .4cm

\pr
The equality is again an application of Theorem 2.13.
The condition for density is necessary, since point evaluations are continuous.
The sufficiency is immediate if one recalls that
$C_c^\infty(U)$ is dense in $C^m(U)$ ([\treref, Theorem 15.3]), showing
that the subspace in 1 is dense if $f_0$ has no zeros.
\qed

\vskip .4cm

\par\noindent
If we take $f_0=1$ and $\Phi$ the identity, then we obtain the well-known
density
of the polynomials in $C^m(U)$ [\treref, Corollary 4, p. 160].
\par\noindent
The Stone-Weierstrass theorem for polynomials also follows from the theorem
(although
this is admittedly not the shortest way to prove it).
Indeed, let $K\subset\R^n$ be compact and non-empty. Choose an open
neighbourhood
$U$ of $K$. If $f\in C_c(K)$, then $f$ has a continuous extension $f^{ext}\in
C(U)$ as a consequence of Tietze's theorem [\bref, Theorem 10.4, p. 30].
Since the polynomials are dense in $C(U)$, $f^{ext}$ can be approximated
uniformly
by polynomials on any compact subset of $U$, in particular on $K$.

\vskip .4cm

\par\noindent
The same proof as for Theorem 3.2 yields:

\vskip .4cm

\par\noindent
{\bf Theorem 3.3.} {\sl Let $f_0\in\schw(\R^n)$ and let $\Phi:\R^n\mapsto
\Phi(\R^n)$ be a diffeomorphism of class $C^\infty$. Suppose that there exists
$\e>0$ such that
$$
\sup_{x\in\R^n} \Vert x\Vert^N \left\vert \ex {\e\Vert\Phi\Vert}
                 g D^\a f_0 \right\vert <\infty.
$$
for all $N\geq 0$, all multi-indices $\a$ and all functions $g$ that are
polynomials in arbitrary derivatives of the components of $\Phi$.
Then the annihilators of the following subspaces of $\schw(\R^n)$ are equal:
\item{1.}$\{\left(f\circ\Phi\right)f_0\st f\in\schw(\R^n)\}$.
\item{2.}$\fialg f_0$.
\item{3.}${\rm Span}\{\ex {i(\l,\Phi)} f_0\st\l\in\C_\e^n \}$.
\item{4.}${\rm Span}\{\ex {i(\l,\Phi)} f_0\st\l\in\R^n \}$.
\par\noindent
These subspaces are dense in $\schw(\R^n)$ if and only if $f_0$ has no zeros.}

\vskip .4cm

\par\noindent
The density of the span of the Hermite functions, as \lq\lq proved" in the
Introduction, follows
from the application of the theorem with $f_0(x)=\exp (-x^2/2)$ and $\Phi$ the
identity.

\vskip .4cm

\par\noindent
Let $G$ be a locally compact abelian group $G$ and $f_0\in L_1(G)$. Then it is
well known that the translates of $f_0$ span a dense subspace of $L_1(G)$
if and only if the Fourier transform of $f_0$ has no zeros
[\rureftwo, Theorem 7.2.5.d].
The following corollary has the same flavour; it shows e.g.\ that the
translates of the Gaussian $\exp( -\vert x\vert^2)$ span a dense subspace
in $\schw(\R^n)$.
The corollary follows from Theorem 3.3 if
one takes $\Phi$ the identity and recalls that the Fourier transform induces a
homeomorphism in $\schw(\R^n)$.

\vskip .4cm

\proclaim Corollary 3.4. Let $f_0\in\schw(\R^n)$. Suppose that there
exists $\e>0$ such that
$$
\sup_{x\in\R^n} \left\vert \ex {\e\Vert x\Vert} D^\a f_0\right\vert<\infty
$$
for all multi-indices $\a$. Then the translates
of the Fourier transform of $f_0$ span a dense subspace of $\schw(\R^n)$ if and
only if
$f_0$ has no zeros.

\vskip .4cm

\xdef\fourfour{\the\pageno}

\newsection {4. Closing remarks}

\vskip .4cm

\par\noindent
The method that we used to exhibit dense subspaces has some flexibility.
Let us take a closer look at the structure of the proof, indicating a possible
method of proof in cases that are not covered by Theorem 2.13.

\vskip .4cm
\par\noindent Let $E$ be some function space on $\R^n$. Fix $f\in E$,
$T\in\edual$,
and consider the map $H_T:\R^n\mapsto\C$ defined by $H_T(\l)=
\langle T\,,\,\ex {-i(\l,\,.\,)} f\rangle$
(we take the identity for $\Phi$ for convenience and assume that
$\ex {-i(\l,\,.\,)} f\in E$ for all $\l\in\R^n$).
Note that the domain of $H_T$ is $\R^n$ and not $\C^n$. What one
{\it really} wants to prove is that $\langle T\,,\,pf\rangle=0$
for all polynomials $p$ is equivalent to $H_T$ being identically zero on $\R^n$
(rather than
$\C^n)$, since $\R^n$ is the domain of the integral in Lemma 2.12.
In our case we proved the (formally obvious) fact that $\langle
T\,,\,pf\rangle$
is in fact a multiple of a derivative of $H_T$,
evaluated at zero --- and $H_T$ happened to extend to a holomorphic map on
$C^n_\e$. But there are other
theorems that ascertain that a function is identically zero if all derivatives
at a point vanish: e.g., if $n=1$ then
one might try to prove that $H_T$ is in a quasi-analytic class (see
[\rurefthree]).
Once this hurdle is taken one can consider the weak integrals as above, try
to prove that they exist and (hopefully) conclude that $T=0$, or at least
obtain
a useful description of $({\rm Span}\{pf\st\hbox{$p$ a polynomial}\})^\perp$.

\vskip .4cm

\par\noindent
As illustration of another way of concluding that $H_T=0$, let us prove the
following proposition. The proof is a variation on Hamburger's method
in [\hamref]. There is a holomorphic function involved, and it follows
immediately
from the hypotheses that many of its derivatives vanish in $0$. As suggested
above,
we use additional information to conclude that the function is in fact equal to
zero.

\vskip .4cm

\proclaim Proposition 4.1. Let $1\leq p<\infty,$ and let $N\geq 0,\,l\geq 2$ be
integers. Then the span of $\{x^n\ex {-x} \st n\geq N,$ $
\,l\nmid n \}$ is dense in $L_p((0,\infty),dx)$.

\vskip .4cm

\pr
Let $g\in L_q((0,\infty),dx)$ (where $q$ is the conjugate exponent of $p$)
and suppose that $g$ is in the annihilator of ${\rm Span}\{x^n\ex {-x} \st
n\geq N,\,l\nmid n \}$.
Put $\O=\{\l\in\C\st {\rm Im}\l<{1\over 2}\}$, and let
$$
H_g(\l)=\int_0^\infty g(x) \ex {-i\l x} \ex {-x} \,dx\quad (\l\in\O).
$$
Then $H_g$ is holomorphic and bounded on $\O$, and
$({d\over d\l})^nH_g (0)=0\,\,(n\geq N,\,\,l\nmid n)$.
Put $\o=\ex {2\p i\over l} $.
There exists a polynomial $P(\l)$ such that $H_g(\o^k\l)-P(\o^k\l)=
H_g(\l)-P(\l)$ for all $\l\in\O$ and all integers $k$ such that $\o^k\l\in\O$
(take the first $N$ terms of the power series of $H_g$ around $0$ for $P$).
This implies that $H_g$ can be extended to an entire function. Indeed, let
$\l\in\C$. Choose $k$ such that $\o^k\l\in\O$ (which is possible since $l\geq
2$)
and put $H_g^{ext}(\l)=H_g(\o^k\l)-P(\o^k\l)+P(\l)$, which is well-defined in
view
of the above. Since $H_g$ is bounded on $\O$, $H_g^{ext}$ is apparently entire
and of polynomial growth, hence in fact equal to a polynomial.
On the other hand,
the restriction of $H_g^{ext}$ to $\R$ is the Fourier transform of $g\ex {-x}
\chi_{(0,\infty)}$, which tends to zero as $\l\rightarrow\pm\infty$ by the
Riemann-Lebesgue lemma. We conclude that
$H_g=0$. But then $g=0$ almost everywhere on $(0,\infty)$ by the injectivity
of the Fourier transform and the fact that $\ex {-x} $ has no zeros.
\qed

\vskip .4cm

\par\noindent
We end with a reformulation of the results in Section 2 in terms of
representation theory.
\par\noindent
The local complex Lie group
$\C_{\e}^n$
acts on $e_0$ {\it as if} the action came from a representation: $\ex
{-i(\l_1,\Phi)}
(\ex {-i(\l_2,\Phi)}  e_0)$$=\ex {-i(\l_1+\l_2,\Phi)} e_0$, provided that
$\l_1,\,
\l_2,\,\l_1+\l_2\in\C_\e^n$ (a consequence of Assumption 2.6.1). But this
action
is in general not the restriction of the obvious candidate for a global action
of $\C_\e^n$ on $E$ since $\ex {-i(\l,\Phi)} e$ simply need not be in $E$ for
all $e\in E$ and all $\l\in\C_\e^n$. One might call $e_0$ a local
representation vector, where the term \lq\lq local" has a double meaning: it
expresses
the fact that $\C_\e^n$ is a local Lie group and also the fact that the
action of this local Lie group is not necessarily globally defined
on $E$. Assumption 2.6 implies more than just this: $e_0$ is in fact
a {\it holomorphic} vector, and the action of the Lie algebra of $\C_\e^n$ on
$e_0$ can be identified. Finally, the density parts
of the theorems in
Section 3 are theorems stating that certain holomorphic
local representation vectors are cyclic for the local action of $\C_\e^n$.

\vskip .4cm

\xdef\fourrefs{\the\pageno}
\newsection {References}
\vskip .2cm

\item{[\abref]}  Aliprantis, C.D., Burkinshaw, O.: {\it Principles of Real
Analysis.}
                 Edward {\hyphenation{Ar-nold}Arnold} Ltd., London, England,
1981.
\item{[\bref]}   Bredon, G.E.: {\it Topology and Geometry.} Springer Verlag,
New York,
                 1993.
\item{[\hamref]} Hamburger, H.: {\it Beitr\"age zur Konvergenztheorie der
                 Stieltjesschen\hfill\break {\hyphenation {Ketten-bru-che}
Kettenbr\"uche}}. Math.\ Z.\ {\bf 4} (1919), 186-222.
\item{[\rurefone]} Rudin, W.: {\it Functional Analysis.} McGraw-Hill, New York,
1991.
\item{[\rureftwo]} Rudin, W.: {\it Fourier Analysis on Groups.} Wiley \& Sons,
New York,
                   1990.
\item{[\rurefthree]} Rudin, W.: {\it Real and Complex Analysis.} McGraw-Hill,
New York,
                     1985.
\item{[\schwref]}  Schwartz, L.: {\it Th\'eorie des distributions.} Hermann,
Paris, 1978.
\item{[\treref]}   Treves, F.: {\it Topological Vector Spaces, Distributions
and Kernels.}
                   Academic Press, New York London, 1967.
\item{[\zaref]}    Zaanen, A.C.: {\it Integration.} North-Holland Publishing
Company,
                   Amsterdam, 1967.

\vfill
\eject
\end